\documentclass[10pt]{article}
\usepackage{amsmath,amsfonts}

\newcommand{\sh}[1]{\!#1\!}

\newcommand{\CCC}{\mathbb{C}}
\newcommand{\kk}{\mathbf{k}}
\newcommand{\ZZ}{\mathbb{Z}}
\newcommand{\PP}{\mathbb{P}}
\newcommand{\NN}{\mathcal{N}}
\newcommand{\MM}{\mathcal{M}}
\newcommand{\LL}{\mathcal{L}}
\newcommand{\Gr}{\mathop{\rm Gr}\nolimits}
\newcommand{\Fl}{\mathop{\rm Fl}\nolimits}
\newcommand{\Span}{\mathop{\rm Span}\nolimits}
\newcommand{\ord}{\mathop{\rm ord}\nolimits}
\newcommand{\rank}{\mathop{\rm rank}\nolimits}
\newcommand{\diag}{\mathop{\rm diag}\nolimits}
\newcommand{\id}{\mathop{\rm id}\nolimits}
\newcommand{\vdim}{\mathop{\rm vdim}\nolimits}
\newcommand{\Lam}{\Lambda}
\newcommand{\Ldot}{\Lambda_{\bullet}}
\newcommand{\Edot}{E_{\bullet}}

\newcommand{\bb}{\mathbf{b}}
\newcommand{\cc}{\mathbf{c}}
\newcommand{\dd}{\mathbf{d}}

\newcommand{\rr}{\mathbf{r}}
\newcommand{\tW}{\widetilde{W}}
\newcommand{\tM}{\tilde{M}}
\newcommand{\ltimes}{{\triangleright\!\!\!\times}}
\newcommand{\at}{{\small $\bigcirc$} \hspace{-1.66em}
 {\tt a}\hspace{.075em}}

\newcommand{\ii}{\mathbf{i}}
\newcommand{\csupset}[1]{\!\stackrel{#1}{\supset}\!}

\newcommand{\GLhat}{\mathop{\hat{\rm GL}}\nolimits}

\newcommand{\GL}{\mathop{\rm GL}\nolimits}

\begin{document}

\centerline{\Large\bf Affine Schubert Varieties
and Circular Complexes}
\vspace{1em}
\centerline{\bf Peter Magyar}
\vspace{.2em}
\centerline{\tt magyar\at\,math.msu.edu}
\vspace{.2em}
\centerline{\tt www.math.msu.edu/\~{}magyar}
\vspace{1em}
\centerline{August 1999}
\vspace{1em}

\noindent 
Schubert varieties have been exhaustively
studied with a plethora of techniques:
Coxeter groups, explicit desingularization, 
Frobenius splitting, etc.
Many authors have applied these techniques
to various other varieties, usually defined by determinantal
equations.  It has turned out that most of these
apparently different varieties are actually Schubert
varieties in disguise, so that one may use
a single unified theory to understand many large families
of spaces. 

The most powerful result in this direction was given by 
Lusztig \cite{lusztig1}, \cite[\S11]{lusztig2},
as a footnote to his work on canonical bases.
He showed that the variety of nilpotent representations
of a cyclic quiver (including nilpotent conjugacy
classes of matrices) is isomorphic to an open subset 
of a Schubert variety for the loop group 
$\widehat{GL}_n$.  In this paper, we 
attempt to describe the affine Schubert varieties (\S1) 
and Lusztig's isomorphism (\S2) in the simplest terms possible.

We then apply this isomorphism to an interesting example,
the variety of circular complexes, recovering many of the
results of Mehta and Trivedi \cite{mehta}. 
(The reader may skip to this application in \S3 immediately
after reading \S1.)
Our technique is similar to that of Lakshmibai and Magyar
\cite{lakshmibai}: it is as a chapter in
the ``ubiquity of Schubert varieties.''

\section{Affine flag variety and Weyl group}

We begin by describing the loop group and its
flag variety as a classical group. 
For more details of
the material of this section see
Pressley-Segal \cite{pressley}, 
Kac-Raina \cite{kac}, 
Slodowy \cite{slodowy}, 
Kazhdan-Lusztig \cite{kazhdan}, 
Kumar \cite[Appendix C]{kumar}, 
Shi \cite{shi}, 
Bjorner-Brenti \cite{bjorner},
and Eriksson-Eriksson \cite{eriksson}.

\subsection{Loop group and affine flag variety}

Let $\kk$ be an arbitrary field, and  
$F:=\kk((t))$, the field of formal Laurent series
$f(t)=\sum_{i\geq N}a_i t^i$ with $a_i\in\kk$; 
and $A:=\kk[[t]]$, the ring of formal Taylor series.
For such $f(t)\sh\neq 0$, we let $\ord(f)$ be the
smallest integer $N$ for which $a_N\neq 0$.

Fix a positive integer $n$, and define 
$G=\GLhat_n(\kk):=\GL_n(F)$,
the group of invertible $n\times n$ matrices with
coefficients in $F$.  
We call this the {\it loop group} because
for $\kk=\CCC$ we may think of $G$ 
as a completion of the group of
polynomial maps from the circle 
$S^1\subset \CCC^\times$ to $\GL_n(\CCC)$.
 
Let $G_j:=\{g\in G\mid \ord\det g=j\}$,
so that $G_j G_k=G_{j+k}$, and for any
$\sigma\in G_1$, we have $G_j=\sigma^j G_0=G_0 \sigma^j$,
and 
$
G=\coprod_{j\in\ZZ} G_j.
$
This should be thought of as the decomposition of $G$ into connected
components.  (For $\kk=\CCC$, and $g$ a polynomial map,
the number $\ord\det g$ is the
winding number of the loop $\det g:S^1\to \CCC^\times$,
and the $G_j$ are the connected components of $G$
in the appropriate compact-open topology.)

Let $V:=F^n$, a vector space over $F$ with a natural
action of $G$. Let $e_1,\ldots, e_n$ denote the standard
$F$-basis of $V$, and for $c\in \ZZ$,
define $e_{i+nc}:=t^c e_i$.  (Thus, $\{e_i\}_{i\in \ZZ}$
is a $\kk$-basis of $V$, in the sense appropriate
to a topological vector space with the $t$-adic topology.)

An {\it $A$-lattice}
$\Lam\subset V$ is the $A$-submodule $\Lam=Av_1\oplus\cdots\oplus
Av_n$, where $\{v_1,\ldots,v_n\}$ is an $F$-basis of $V$.
We may write $\Lam=\Span_\kk\langle v_i\rangle_{i\geq 1}$,
the space of infinite $\kk$-linear
combinations of the vectors $v_{i+nc}:=t^c v_i$.
Consider the family of {\it standard $A$-lattices}:
$$
E_j:=\Span_A\langle e_j, e_{j+1},\ldots
e_{j+n-1}\rangle=\Span_\kk\langle e_i\rangle_{i\geq j}.
$$
Note that $E_j=\sigma^j E_1$, where 
we use the shift operator
$\sigma(e_i):=e_{i+1}$, or as a matrix:
$$
\sigma=
\left(\begin{array}{@{\!}c@{\!}c@{\!}c@{\!}c@{\!}c@{\!}}
\ 0&\ \ 0&\ \cdots&\ 0&\ \ t\ \\[-.3em]
\ 1&\ \ 0&\ \cdots&\ 0&\ \ 0\ \\[-.3em]
\ 0&\ \ 1&\ \cdots&\ 0&\ \ 0\ \\[-.6em]
\ \cdot&\ \ \cdot&&\ \cdot&\ \ \cdot\ \\[-.8em]
\ \cdot&\ \ \cdot&&\ \cdot&\ \ \cdot\ \\[-.4em]
\ 0&\ \ 0&\ \cdots&\ 1&\ \ 0\ 
\end{array}\right)\ 
\in\ G_1
$$

The {\it affine Grassmannian} $\Gr(V)$ is the space
of all $A$-lattices of $V$.  Clearly $\Gr(V)$ is a
homogeneous space with respect to the obvious action
of $G$, and the stabilizer of the standard lattice $E_1$
is $P_{\hat 0}:=GL_n(A)$, 
the subgroup of matrices with coefficients in $A$ 
and with determinant having $\ord=0$.
Thus $\Gr(V)\cong G/P_{\hat 0}$, and
the connected components of the Grassmannian are
$\Gr_j(V):=G_0\cdot E_j=G_j\cdot E_1\cong G_j/P_{\hat 0}$.
In fact, $\Gr_j(V):=\{\Lam\mid \vdim(\Lam)=j\}$,
where we define the {\it virtual dimension }
$$\vdim(\Lam):=\dim_\kk(\Lam/\Lam\cap E_1)
-\dim_\kk(E_1/E_1\cap \Lam) .$$

The {\it complete affine flag variety}
$\Fl(V)$ is the space of all flags of lattices
$\Ldot=(\Lam_1\supset\cdots\supset\Lam_n)$
such that $\Lam_n\supset t\Lam_1$ and
$\dim_\kk(\Lam_j/\Lam_{j\sh+1})=1$.
There always exists an $F$-basis $\{v_1,\ldots,v_n\}$ of $V$
such that 
$\Lam_j = \Span_A\langle v_j,\ldots,v_n,
tv_1,\ldots,tv_{j\sh-1}\rangle$
$=\Span_\kk\langle v_i\rangle_{i\geq j}$,
where $i$ runs over all integers not less than $j$,
and $v_{i+nc}:=t^c v_i$.
The {\it standard flag} is $\Edot:=(E_1\supset\cdots\supset E_n)$,
whose stabilizer $B$ is the subgroup of 
matrices $b\in P_{\hat 0}$ which are lower-triangular modulo
$t$:
$$
B:=\{b=(b_{ij})\in GL_n(A)\mid \ord(b_{ij})\sh>0\ \
\forall\, i\sh<j\}.
$$
Thus, $\Fl(V)\cong G/B$, with connected components
$\Fl_j(V):=G_j\sh\cdot\Edot\cong G_j/B=
\{\Ldot\mid \vdim(\Lam_1)=j\}$.
Furthermore, the projection $\Fl(V)\to\Gr(V)$, 
$\Ldot\mapsto\Lam_1$ is a bundle whose fiber
is the space of complete flags in the $n$-dimensional
$\kk$-vector space $\Lam_1/t\Lam_1$.

\subsection{Affine Weyl group}

We first discuss $\tW$, the Weyl group of the disconnected
group $G$; and then $W$,
the Weyl group  of the connected component $G_0$.

Let $S_\infty$ be the group of bijections $\pi:\ZZ\to\ZZ$,
and let $\sigma:i\mapsto i\sh+1$ be the shift
bijection.  
Let $\tW\subset S_\infty$ be the subgroup
of bijections which commute with the $n$th power of $\sigma$: 
that is, $\tW:=\{\pi\in S_\infty\mid \pi\tau=\tau\pi\}$,
where $\tau:=\sigma^n:i\mapsto i\sh+n$.

For $\cc=(c_1,\ldots,c_n)\in\ZZ^n$, define an element
of $\tW$,\  $\tau^\cc:
i\mapsto i\sh+nc_{(i\!\!\!\mod n)}$.  
This gives an
embedding of the additive group $\ZZ^n\subset \tW$.
Furthermore, we have the embedding $S_n\subset \tW$
with $\bar\pi(i+nc):=\bar\pi(i)+nc$ for $\bar\pi\in S_n$.
Then we may write any element
$\pi\in\tW$ as $\pi=\bar\pi \tau^\cc$
for unique $\bar\pi\in S_n$, $\cc\in\ZZ^n$,
and we have 
$\bar\pi_1\tau^{\cc_1}\bar\pi_2\tau^{\cc_2}
= \bar\pi_1\bar\pi_2
\tau^{\bar\pi_2^{-1}(\cc_1)+\cc_2}$.
That is, $\tW= S_n\ltimes \ZZ^n$, a
semi-direct product.  The normal subgroup
$\ZZ^n$ is the kernel of the homomorphism
$\tW\to S_n$ which takes each $\pi:\ZZ\to\ZZ$
to a permutation of cosets $\bar\pi:\ZZ/n\ZZ\to\ZZ/n\ZZ$.

Thus, an element $\pi\in \tW$
is equivalent to a sequence of integers
$[\pi(1),\ldots,\pi(n)]$ such that $i\mapsto \bar\pi(i)$
defines a permutation of $\ZZ/n\ZZ$.
For example, in this one-line notation 
$\sigma=[2,3,\cdots,n\sh+1]$ and 
$\tau^\cc=[1+nc_1,2+nc_2,\ldots,n+nc_n]$.
  
We obtain another useful notation by embedding $\tW\subset G$.
If we let $\pi=\bar\pi\tau^\cc\in W$ act $F$-linearly on $V$ by 
$\pi(e_i):= e_{\pi(i)}$, the corresponding
matrix is the {\it affine permutation matrix} $(a_{ij})$
with $a_{\bar\pi(i),i}=t^{c_i}$.
For example, $\sigma$ becomes the matrix in $G_1$ 
of the previous section,  
$\tau^\cc=\diag(t^{c_1},\cdots,t^{c_n})$,
and $\tau=\tau^{(1,\cdots,1)}=\diag(t,\cdots,t)$,

Using this embedding, we may show that
$\tW\cong N_G(\hat T)/\hat T$, 
where $\hat T$ is the subgroup of diagonal matrices 
with entries in $A$.

Considering the embedding $\tW\subset G$, we define
$$
W_j:=\tW\cap G_j
=\left\{\pi\in\tW\ \left|\ 
\sum_{i=1}^n\pi(i)-i=j\right.\right\},\qquad
W:=W_0.
$$
Thus $W$ is a normal subgroup of $\tW$, and 
we have $W= S_n\ltimes\ZZ^n_0$,
where $\ZZ^n_0:=\ZZ^n\cap W=\{\cc\mid \sum_{i=1}^n c_i=0\}$.
That is, $W\cong \widehat{S}_n$, the affine Weyl group
of extended Dynkin type $\widehat A_{n\sh-1}$.

Define the simple reflections $s_0,\ldots,s_{n\sh-1}$ in $W$ 
by $s_i(i)=i\sh+1$, $s_i(i\sh+1)=i$, and $s_i(j)=j$ for 
$j\not\equiv i,\, i\sh+1\!\!\mod n$.  
We shall sometimes denote $s_n:=s_0$.
Then $s_0,\ldots,s_{n\sh-1}$ are involutions
generating $W$ and satisfying the
Coxeter relations $(s_i s_{i\sh+1})^3
=\id$ for $0\leq i\leq {n\sh-1}$, and $(s_i s_j)^2=\id$ otherwise.
We have a semi-direct
product $\tW=\langle\sigma\rangle \,\ltimes W$.  
Here $\sigma$ acts on $W$ via the outer automorphism:
$\sigma s_i\sigma^{-1}=s_{i\sh+1}$.

The Bruhat length $\ell(\pi)$ 
is defined as usual for $\pi\in W$, as the
smallest number of simple reflections whose product is
$\pi$;
and we extend this to $\tW$ by letting
$\ell(\sigma^j\pi):=\ell(\pi)$.
We have J.~Shi's formula \cite{shi}: 
$$\ell(\pi)=\sum_{1\leq i<j\leq n} 
\left|\,\mathop{\rm floor}(
\tfrac{\pi(j)-\pi(i)}{n})\right|,
$$
where $\mathop{\rm floor}(x)$ denotes the greatest integer
not exceeding $x$.
Lusztig showed:
$$
\ell(\tau^\cc)=(n\sh-1)c_1+(n\sh-3)c_2+\cdots
+(-n\sh+3)c_{n-1}+(-n\sh+1)c_n\, ,
$$
namely, the dot product of $\cc$ with
$$2\rho^{\vee}:=(n\sh-1,n\sh-3,\ldots,-n\sh+3,-n\sh+1)
=\!\!\sum_{1\leq i<j\leq n}\!\!\!\! e_i-e_j.$$

\subsection{Wiring diagrams}

The structure of the Weyl group is further elucidated by
the {\it loop wiring diagrams} (cf.~Berenstein-Fomin-Zelevinsky
\cite{bfz}).  Consider
a cylinder $[0,1]\times S^1$.  On the right end,
label each point $(1,e^{2\pi\sqrt{-1}i/n})$ 
with the integer $i$,
and similarly on the left end.  
Now, represent a permutation
$\pi=\bar\pi\tau^\cc\in \tW$ 
by $n$ curves, each joining a
point $i$ on the right to the point $\bar\pi(i)$
on the left, but looping counter-clockwise 
around the cylinder $c_i$ times.   
\\[1em]
{\bf Example.}  For $n=3$, the permutation
$\pi=[-2,2,6]\in W$, with $\bar\pi=\id\in S_3$ and
$\cc=(-1,0,1)\in\ZZ^3_0$, 
is represented by the picture:

\bigskip

%TexCad Options
%\grade{\off}
%\emlines{\off}
%\beziermacro{\on}
%\reduce{\on}
%\snapping{\off}
%\quality{2.00}
%\graddiff{0.01}
%\snapasp{1}
%\zoom{1.00}
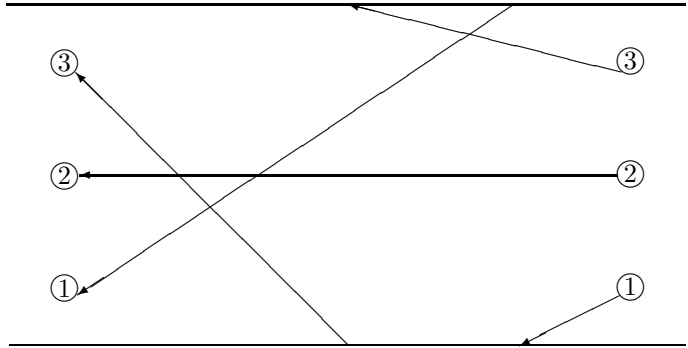
\unitlength .75mm
\linethickness{0.4pt}
\begin{picture}(130.00,70.42)
\put(130.00,70.33){\line(-1,0){120.33}}
\put(120.33,60.33){\circle{4.47}}
\put(20.00,60.00){\circle{4.71}}
\put(120.33,40.33){\circle{4.47}}
\put(20.00,40.00){\circle{4.71}}
\put(120.33,20.33){\circle{4.47}}
\put(20.00,20.00){\circle{4.71}}
\put(130.00,10.00){\line(-1,0){119.67}}
\put(120.33,20.33){\makebox(0,0)[cc]{1}}
\put(20.00,20.00){\makebox(0,0)[cc]{1}}
\put(20.00,40.00){\makebox(0,0)[cc]{2}}
\put(120.33,40.33){\makebox(0,0)[cc]{2}}
\put(20.00,60.00){\makebox(0,0)[cc]{3}}
\put(120.33,60.33){\makebox(0,0)[cc]{3}}
\put(99.67,70.33){\vector(-3,-2){77.33}}
\put(118.67,58.33){\vector(-4,1){48.33}}
\put(70.33,10.00){\vector(-1,1){48.33}}
\put(118.00,40.00){\vector(-1,0){95.33}}
\put(118.33,18.67){\vector(-2,-1){17.67}}
\end{picture}

\vspace{-1em}

\noindent Here we represent our cylinder by identifying the
top and bottom borders of the picture, so that each point
$i$ on the right is connected to the same $i$ on
the left (since $\bar\pi=\id$).
However, since $\cc=(-1,0,1)$, 
the curve starting from 1 travels once
clockwise around the cylinder, the curve from 2
travels straight across, and the curve from 3 travels once
counter-clockwise.

We may read off much combinatorial data from this picture.
Since the curves have a total of 4 crossings, we conclude
that $\ell(\pi)=4$. By listing these crossings, 
as well as crossings over the top and bottom margins, 
we obtain a reduced decomposition: 
$\pi=s_2s_1s_2\sigma s_2\sigma^{-1}$.  
That is, the leftmost crossing switches the top two curves, 
giving a factor $s_2$; the second switches the bottom 
two curves, $s_1$; again $s_2$; then the bottom curve 
crosses to the top, $\sigma$; 
again $s_2$; and finally the top curve
crosses to the bottom, $\sigma^{-1}$.
Using $s_i \sigma=\sigma s_{i\sh-1}$,
we have $\pi=s_2s_1s_2s_0$.

\subsection{Schubert varieties}

By Gaussian elimination, we obtain the 
{\it Bruhat decomposition} of $G$ into double $B$-cosets:
$G=\coprod_{\pi\in \tW} B\pi B$, where we consider each
$\pi$ as an affine permutation matrix.
Hence we also have a Bruhat decomposition of
the affine flag variety
$\Fl(V)=\coprod_{\pi\in\tW} X^\circ_\pi$
into {\it Schubert cells} $X^\circ_\pi:= B\sh\cdot \pi \Edot$,
where $\pi \Edot$ is a translation of the standard flag
$\Edot=(E_1\supset\cdots\supset E_n)$.
In particular, $\Fl_j(V)$ is the union of 
all $X^\circ_\pi$ with $\pi\in W_j$.  The Schubert
cells can be defined by dimension
constraints called {\it Schubert conditions}.
For $\pi\in W_j$, we have:
\begin{eqnarray*}
X^\circ_\pi&=&
\left\{\Ldot\in\Fl(V)\left|\,\begin{array}{c}
%\dim_\kk(\Lam_i/\Lam_i\sh\cap E_j)
%=\#(\pi\ZZ_{\geq i}\sh\setminus \ZZ_{\geq j})\\ 
\dim_\kk(E_j/\Lam_i\sh\cap E_j)
=\#(\ZZ_{\geq j}\sh\setminus \pi\ZZ_{\geq i})
\end{array}\right.\right\}
\\[.3em]
&=&
\left\{\Ldot\in\Fl_j(V)\left|\,\begin{array}{c}
\dim_\kk(\Lam_i/\Lam_i\sh\cap E_j)
=\#(\pi\ZZ_{\geq i}\sh\setminus \ZZ_{\geq j}) 
\end{array}\right.\right\}
\end{eqnarray*}
where $\ZZ_{\geq i}$ denotes the integers not less than
$i$; we define $\pi\ZZ_{\geq i}:=
\{\pi(i),\pi(i\sh+1),\ldots\}$;
and $\setminus$ denotes set complement.
Indeed, the set on the right of the equation
is clearly $B$-invariant,
and $\pi'\Edot$ lies in this set 
if and only if $\pi'=\pi$.  

The {\it Schubert variety}, meaning the
topological closure $X_\pi :=\overline{X^\circ_\pi}$, is
obtained by replacing $=$ in the above 
Schubert conditions with $\leq$.
We say $\pi\leq\pi'$ in the {\it Chevalley-Bruhat order} 
if $X_\pi\subset X_{\pi'}$, and
we can express this combinatorially as:
$\pi\leq \pi'$ iff \
$\#\pi\ZZ_{\geq i}\sh\setminus \ZZ_{\geq j}\leq
\#\pi'\ZZ_{\geq i}\sh\setminus \ZZ_{\geq j}$\
for all $1\leq i\leq n$,\ $j\in\ZZ$.

We explain below how $X_\pi$ has the structure
of a projective algebraic variety.  
With this structure, 
the varieties $X_\pi$ include as special cases
the familiar Schubert varieties for $\GL_n(\kk)$.
In fact, for $\pi=\bar\pi\in S_n\subset\tW$
and $\Ldot\in X_{\pi}$, we have 
$\pi(\ZZ_{\geq 1})=\ZZ_{\geq 1}$ 
and $\pi(\ZZ_{\geq n\sh+1})=\ZZ_{\geq n\sh+1}$.
Hence $\Lam_1/\Lam_1\sh\cap E_1=0$ and $\Lam_1\subset E_1$.  
Also $n=\dim_\kk(\Lam_1/\Lam_1\sh\cap E_{n\sh+1})\leq
\dim_\kk(E_1/E_{n\sh+1})=n$, so $\Lam_1=E_1$ and
$\Lam_i\supset t\Lam_1=E_{n\sh+1}$ for $1\leq i\leq n$.
Letting $\kk^n=E_1/E_{n\sh+1}$ and $V_i:=\Lam_i/E_{n\sh+1}$, 
we thus find that $\Ldot\in X_\pi$ is in natural correspondence
with the complete flag  
$$
\kk^n=V_1\supset V_2\supset\cdots\supset V_n\supset 0,
$$
and the affine Schubert conditions on $\Ldot$
are equivalent to the usual Schubert conditions
$$
\dim_\kk(V_i\cap \bar E_j)\geq\#(\pi[i,n]\cap [j,n])
$$
relative to the standard flag $\bar E_j:=E_j/E_{n\sh+1}$.
For example, $X_{\id}=\{\Edot\}$, a single point.

For a general $\pi\in \tW$, we can find $a, b$ so that
$\ZZ_{\geq a}\supset \pi\ZZ_{\geq i}\supset \ZZ_{\geq b}$ 
for $1\leq i\leq n$.  Then any $\Ldot\in X_\pi$ 
satisfies $E_a\supset\Lam_i\supset E_b$, and we 
may embed $X_\pi$ inside a partial flag variety of
the finite-dimensional $\kk$-vector space $E_a/E_b$.
The flags in the image of this embedding must
satisfy certain ordinary Schubert conditions,
but they must also be stable under the nilpotent
map induced on $E_a/E_b$ by $t\in A$.
This makes $X_\pi$ into an algebraic variety
over $\kk$ (in fact, even defined over the integers).

We can imitate all the standard geometric constructions
for Schubert varieties of $\GL_n(\kk)$.
For example, we can show $\dim_\kk X_\pi=\ell(\pi)$; 
we can explicitly construct Bott-Samelson
resolutions of $X_\pi$ as configuration varieties; 
and we can use the usual Frobenius-splitting 
arguments to show that the variety 
$X_\pi$ is normal, Cohen-Macaulay, etc.
\\[1em]
{\bf Example.} Consider as above $n=3$,
$\pi=[-2,2,6]$.  Take $a=-2$, $b=4$,
and write: 
$$
X_\pi=\{\ (E_{-2}\stackrel{3}{\supset}
\Lam_1\stackrel{1}{\supset}
\Lam_2\stackrel{1}{\supset}
\Lam_3\stackrel{1}{\supset} E_4)\ 
\left|\begin{array}{c}
t\Lam_1\subset\Lam_3 ,\
\Lam_2\subset E_1
\end{array}\right.\},
$$
\vspace{-1em}

\noindent
where $U\stackrel{d}{\supset} V$ means $U\supset V$
and $\dim_\kk(U/V)=d$.
All the Schubert conditions for $X_\pi$
follow from the conditions specified on the right side of
the equation.

Let $\kk^6=E_{-2}/E_{4}$ with basis 
$\{\bar e_{-2},\bar e_{-1},\bar e_{0},
\bar e_1,\bar e_2,\bar e_3\}$, and take $\bar t:\kk^6\to\kk^6$, 
$\bar e_i\mapsto \bar e_{i\sh+3}\!\!\mod E_{4}$, so that
$\bar t^{\,2}=0$.  
Then we have the isomorphism
$$
X_\pi\cong\{\
(\kk^6\stackrel{3}{\supset}
V_1\stackrel{1}{\supset}
V_2\stackrel{1}{\supset}
V_3\stackrel{1}{\supset} 0)\ 
\left|\begin{array}{c}
\bar t(V_1)\subset V_3 ,\ 
V_2\subset \bar E_1
\end{array}\right.\},
$$
a subvariety of a partial flag variety of $\GL_6(\kk)$
defined by Schubert conditions
%$\dim(V_1/V_1\sh\cap \bar E_1)\leq 1$ (or equivalently
%$\dim(V_1\sh\cap \bar E_1)\geq 2$) 
and the algebraic 
incidence condition $\bar t(V_1) \subset V_3$.

Further, we can construct a Bott-Samelson 
variety corresponding to the reduced word
$\pi=s_2s_1s_2s_0$:
$$
Z_{2120}:=
\left\{\begin{array}{c}
(\Lam_1,\Lam_2,\Lam_3,\Lam'_3)\\
\in\Gr(V)^4
\end{array}\left|
\begin{array}{c@{\!}c@{\!}c@{\!}c@{\!}c@{\!}c@{\!}c}
&&E_2&\,\leftarrow\,& E_3&&\\
&\swarrow&&\nwarrow&&\nwarrow&\\
E_1&&&&\Lam'_3&\,\leftarrow&\ tE_1\\
&\nwarrow&&\swarrow&&\swarrow&\\
\Lam_1&\,\leftarrow\,&\Lam_2&\,\leftarrow\,&\Lam_3&\,\leftarrow&\ t\Lam_1
\end{array}\right.\right\},
$$
where each arrow $U\sh\leftarrow V$
indicates the condition $U\stackrel{1}{\supset} V$.
We may build up this variety by starting with a
single point (the standard flag) and successively adding the
spaces $\Lam'_3$, $\Lam_2$, $\Lam_3$, $\Lam_1$, corresponding
to the reflections $s_2,s_1,s_2,s_0$.
Clearly $Z_{2120}$ is an iterated $\PP^1$-fibration
(thus smooth), and it maps birationally to $X_\pi$
by dropping $\Lam_3'$.  (In this case,
$X_\pi$ happens to be smooth itself.)  
Note that the pattern of inclusions defining $Z_{2120}$ 
is the dual graph of the wiring diagram (turned sideways).

\subsection{Partial flag variety and opposite cell}

A subset $I\subset [0,n\sh-1]$ corresponds to a
{\it parabolic subgroup} $G\supset P_I\supset B$ with
Weyl group $W_I:=\langle s_i\rangle_{i\in I}\subset W$.  
For $I\ni 0$, the partial flag variety corresponding 
to the complement $\widehat I:=[0,n\sh-1]\sh\setminus I$ is: 
$$
G/P_{\widehat I}\cong
\left\{(\Lam_1\supset\cdots\supset\Lam_h)
\left|\begin{array}{c} 
\Lam_j\ \mbox{\rm a lattice},\ \ \Lam_h\supset t\Lam_1,
\\[.2em]
\dim(\Lam_j\sh\setminus\Lam_{j\sh+1})=i_{j\sh+1}-i_j
\end{array}\right.\right\}\ .
$$

We shall find it convenient to index parabolics
by {\it compositions of $n$}: that is, sequences
of positive integers $\dd=(d_1,\ldots,d_h)$ 
with $d_1+\cdots+d_h=n$.
Given $I=\{0\sh=i_1<\cdots<i_{h}\}$, let $i_{h\sh+1}:=n$, and
define a composition by $d_j:=i_{j\sh+1}-i_j$
(so that $i_j=d_1+\cdots+d_{j\sh-1}$).
Then we may rewrite the above more concisely as:
$$
G/P_{\widehat I}\cong
\Fl(\dd;V):= \{(\Lam_1\stackrel{d_1}{\supset}
\cdots\stackrel{d_{h\sh-1}}{\supset}\Lam_h
\stackrel{d_{h}}{\supset}t\Lam_1) \}\ .
$$
Denoting $W_\dd:=W_{\widehat I}$ and
$E_{(j)}:=E_{1+i_j}$,\
$E_{(\bullet)}:=(E_{(1)}\supset
\cdots\supset E_{(h)})$, we have:
$$
\Fl(\dd;V)=\coprod_{\pi W_\dd \in\tW/W_\dd}
B\sh\cdot \pi E_{(\bullet)}\ .
$$
In particular, for $I=\{0\}$,\ $\dd=(n)$,
we have $\tW/W_\dd=\tW/S_n\cong \ZZ^n$,
and 
$$
\Gr(V)=G/P_{\hat 0}=\coprod_{\cc\in\ZZ^n} X_\cc^\circ\ ,
$$
where
$
X_\cc^\circ:=B\sh\cdot \tau^\cc E_1 =
\{\Lam\mid \dim(\Lam/\Lam\sh\cap E_j)=
\#(\tau^\cc\ZZ_{\geq 1}\sh\setminus \ZZ_{\geq j})\,\}\ .
$

\medskip

Next, for any Schubert variety $X_\pi$,
we define a certain affine open subset, 
the {\it opposite cell} $X'_\pi\subset X_\pi$  
(meaning the opposite to the cell 
$X^\circ_\pi\subset X_\pi$, though $X'_\pi$ 
itself is generally not a topological cell).
Let $E'_k:=\Span_{\kk}\langle e_i\rangle_{i<k}$ be the
complementary space to $E_{k}$.  Note that 
$E'_k$ is {\it not} an $A$-lattice in $V$:
rather, it is a lattice over the ring
$A'=k[t^{-1}]\subset F$. For
$\pi\in W$, define $X'_\pi\subset X_\pi\subset 
\Fl_0(V)$ as the set of $\Ldot\in X_\pi$ such that
$\Lam_i\cap E'_{i}=0$ for $1\leq i\leq n$. 
For example, $\Edot\in X'_\pi$ for any $\pi\in W$.

[Note: The condition $\Lam_i\cap E'_{i}=0$
is equivalent to $\Lam_i\oplus E'_{i}=V$.
Proof:  Recall that 
$\dim_\kk(\Lam_i/\Lam_i\sh\cap E_i)-\dim_\kk(E_i/E_i\sh\cap\Lam_i)=0$,
and let $\phi:\Lam_i/\Lam_i\sh\cap E_i\subset
E'_{i}\sh\oplus E_i/\Lam_i\sh\cap E_i\to
E_i/\Lam_i\sh\cap E_i$.
Thus $\Lam_i\cap E'_{i}=\mathop{\rm Ker}(\phi)=0$
$\Longleftrightarrow$ $\mathop{\rm Im}(\phi)=E_i/\Lam_i\sh\cap E_i$
$\Longleftrightarrow$ $E'_{i}+\Lam_i=E'_{i}+E_i=V$.]

More generally, for $\pi\in W_k=\sigma^k W$,\
$X_\pi\subset \Fl_k(V)$, we let
$$\begin{array}{rcl}
X'_\pi&:=&\{ \Ldot\in X_\pi\mid \Lam_i\cap E'_{i\sh+k}=0,\ \
1\sh\leq i\sh\leq n\}\\[.2em]
&=&\{ \Ldot\in X_\pi\mid \Lam_i\oplus E'_{i\sh+k}=V,\ \
1\sh\leq i\sh\leq n\}.
\end{array}$$
Thus $\sigma^k\Edot\in X'_\pi$.
We define $X'_\pi\subset X_\pi\subset\Fl_k(\dd,V)$
similarly: e.g., for $X'_\pi\subset X_\pi\subset\Gr(V)$,
we require $\Lam\cap E'_{k\sh+1}=0$, so $E_{k\sh+1}\in X'_\pi$.
\medskip

Now we examine certain affine
Grassmannian Schubert varieties
which will occur in the following section. 
Suppose $\cc=(c_1,\ldots,c_n)$ satisfies
$0\leq c_1\leq\cdots\leq c_n\leq n$ and $c_1+\cdots+c_n=n$.
Let $c'_j:=\#(\tau^\cc\ZZ_{\geq 1}\sh\setminus
\tau^j\ZZ_{\geq 1})=\#\{i\mid c_i\leq j\}$, the 
conjugate-complement partition of $\cc$.
Then we have:
$$
X^\circ_\cc=\left\{\Lam\in\Gr(V)\mid
E_1\supset\Lam\supset t^n E_1,\
\dim(\Lam/\Lam\cap t^j\!E_1)=c'_j,\ 1\leq j\leq n
\right\}
$$
Since the maximal parabolic $P_{\hat 0}$ stabilizes $t^j E_1$, 
we have $X^\circ_\cc = P_{\hat 0}\sh\cdot \tau^\cc E_1$. 
Furthermore, letting $X^{\circ\prime}_\cc:=
X^\circ_\cc\cap X'_\cc$, we obtain:
$$
X^{\circ\prime}_\cc=\left\{\Lam\in\Gr(V)
\left|\begin{array}{c}
\dim(\Lam/\Lam\cap t^j\!E_1)=c'_j,\ 1\leq j\leq n
\\[.2em]
E_1\supset\Lam\supset t^n E_1,\quad 
\Lam\cap t^n E'_1=0
\end{array}\right.\right\}.
$$
Since $\GL_n(\kk)\subset P_{\hat 0}$ is the joint
stabilizer of $t^j E_1$ and $t^j E'_1$, we have
$X^{\circ\prime}_\cc=
\GL_n(\kk)\sh\cdot \tau^\cc E_1$.

\pagebreak

\section{Lusztig's isomorphism}

In this and the following section, we consider how
certain varieties of matrices may be considered
as opposite cells in affine Schubert varieties. 

\subsection{Nilpotent matrices}

Let $\NN\subset M_{n\sh\times n}(\kk)$ be the
set of nilpotent $n\sh\times n$ complex matrices,
on which $GL_n(\kk)$ acts by conjugation.
Lusztig \cite{lusztig1} has given an
equivariant algebraic isomorphism between $\NN$ and 
the opposite cell of a Schubert variety in $\Gr(V)$.

A matrix in $GL_n(\kk)$ has a natural $A$-linear
action on $V$, and for $N\in\NN$
we can define 
$\phi_N:V\to V$,\, 
$$
\begin{array}{rcl}
\phi_N(v)&:=&
\displaystyle\frac{t^{n\sh-1}}{1-t^{-1}N}(v)\\[1em]
&=& t^{n\sh-1}v+t^{n\sh-2}N(v)+
t^{n\sh-3}N^2(v)+\cdots
+N^{n\sh-1}(v).
\end{array}
$$
Lusztig's isomorphism is given by the map 
$$
\begin{array}{cccc}
\Phi:&\NN&\to&\Gr(V)\\[.1em]
&N&\mapsto&\phi_N(E_1).
\end{array}
$$
Note that $\Phi$ is $GL_n(\kk)$-equivariant:
for $g\in GL_n(\kk)$, we have $\Phi(gNg^{-1})=g \,
\phi_N(g^{-1}E_1)=g\,\phi_N(E_1)= g\,\Phi(N)$.
We also have $E_1\supset\Phi(N)\supset t^n E_1$
for all $N\in\NN$.

We may parametrize the $GL_n(\kk)$-orbits in $\NN$
by $n$-tuples of integers 
$\bb=(b_1,\ldots,b_n)$,
where $n\geq b_1\geq \cdots\geq b_n\geq 0$
and $b_1+\cdots+b_n=n$. 
That is, the orbit $\NN_\bb\subset\NN$ 
consists of those nilpotents whose largest Jordan
block has size $b_1$, the next largest has size $b_2$, etc.  
The open orbit of principal nilpotents is $\NN_{(n,0,\ldots,0)}$, 
the closed orbit $\{0\}=\NN_{(1,\ldots,1)}$.

Let $\cc=(n-b_1,\ldots,n-b_n)$.
Applying elementary linear algebra
to the description of $X^{\circ\prime}_\cc$ in the
previous section, we may easily show that:
$$
\Phi(\NN_\bb)=X^{\circ\prime}_\cc
\qquad\mbox{and}\qquad
\Phi(\overline{\NN}_\bb)=X'_\cc,
$$
where $\overline{\NN}_\bb$ denotes the closure.
In particular, 
$$
\Phi(\NN)=X'_{(0,n,\cdots,n)}=
\{\Lam\in\Gr(V)\mid \Lam\stackrel{n}{\supset}
t^nE_1\}\, .
$$

\noindent 
{\bf Example.} 
Note that we can renormalize our map by a $\tau$-shift:
$\Phi(\overline{\NN}_\bb)=X'_\cc\cong X'_{\tau^j\cc}$
(equivariant isomorphism) for any $j\in\ZZ$.
Thus for the example of \S3, \S4, we have:
$X'_{(-1,0,1)}\cong X'_{(1,2,3)}\cong 
\Phi(\overline{\NN}_{(2,1,0)})$.
\\[1em]
Let us write our map in coordinates. 
We represent $v=\sum_{i=1}^\infty a_i e_i\in E_1$ 
(where $a_i\in\kk$) 
by the semi-infinite column vector with entries $a_i$; and 
$\Lam=\Span_A\langle v_1,\ldots,v_n\rangle\in\Gr(V)$ 
by the semi-infinite matrix $[v_1,\cdots,v_n]$.
Then we may write
$$
\Phi(N)\quad=\quad
\left\lceil\,\begin{array}{@{\!}c@{\!}}
N^{n\sh-1}\\[-.7em] \vdots\\N^2\\N\\I\\0\\[-.7em] \vdots
\end{array}\,\right\rceil,
$$
where $I$ is the identity matrix.
From this we see how the {\it Plucker coordinates}
(the $n\sh\times n$ minors of this matrix) restrict
to polynomial functions on $\NN$.  In particular,
since the vanishing ideal of
the Schubert subvarieties $X_\cc\subset\Gr(V)$
is generated by the vanishing of certain Plucker coordinates,
we obtain generators for the ideal of
$\overline{\NN}_\bb\subset \NN$.
(Cf.~Weyman \cite{weyman}).

\subsection{Cyclic quivers}

Lusztig \cite{lusztig2} has generalized the 
above isomorphism (and simultaneously another 
isomorphism of Zelevinsky \cite{zelevinsky},
\cite{lakshmibai}). 
The generalization involves a positive integer parameter $h$,
with the case $h=1$ reducing to our discussion
of nilpotent matrices.

The {\it cyclic quiver} $\widehat A_{h\sh-1}$
is the oriented graph:\\[-1.2em]

\hspace{1in}
%TexCad Options
%\grade{\off}
%\emlines{\off}
%\beziermacro{\on}
%\reduce{\on}
%\snapping{\off}
%\quality{2.00}
%\graddiff{0.01}
%\snapasp{1}
%\zoom{1.00}
\unitlength .75mm
\linethickness{0.4pt}
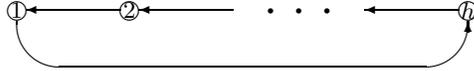
\begin{picture}(91.67,26.00)
\put(10.00,20.00){\circle{3.33}}
\put(30.00,20.00){\circle{3.33}}
\put(90.00,20.00){\circle{3.33}}
\put(28.33,20.00){\vector(-1,0){16.67}}
\put(48.33,20.00){\vector(-1,0){16.67}}
\put(88.33,20.00){\vector(-1,0){16.67}}
\put(30.00,20.00){\makebox(0,0)[cc]{{\small 2}}}
\put(90.00,20.00){\makebox(0,0)[cc]{{\small $h$}}}
\put(65.00,20.00){\circle*{0.94}}
\put(60.00,20.00){\circle*{0.94}}
\put(55.00,20.00){\circle*{0.94}}
\put(10.00,20.00){\makebox(0,0)[cc]{{\small 1}}}
\put(50.00,18.00){\oval(80.00,16.00)[b]}
\put(90.00,16.00){\vector(0,1){2.00}}
\end{picture}
\vspace{-1.5em}

\noindent
For a fixed $h$-tuple of positive integers
$\dd=(d_1,\cdots,d_h)$, 
we define the {\it $\dd$-dimensional representations} of this
quiver to be the affine space 
$$
M_{d(h)\times d(1)}(\kk)\times
M_{d(1)\times d(2)}(\kk)\times\cdots\times
M_{d(h\sh-1)\times d(h)}(\kk).
$$
(For legibility, we have written $d(j)$ instead of $d_j$.)
That is, a representation $(M_1,\ldots,M_h)$ is a way
of replacing each arrow $i\to i\sh-1$
by a linear map $M_i:\kk^{d(i)}\to\kk^{d(i\sh-1)}$,
where we take $d_0:=d_h$.  
(For all $j,k$, we write $d_{j+hk}:=d_j$.)
We have a natural action of the group
$\GL_\dd(\kk):=\GL_{d(1)}(\kk)\times\cdots\times
\GL_{d(h)}(\kk)$ on the space of representations:
$$
(g_1,\ldots,g_h)\cdot (M_1,\ldots,M_h):=
(g_hM_1g_1^{-1},g_1M_2g_2^{-1},\ldots,g_{h\sh-1}M_hg_h^{-1}).
$$

Our main concern is a certain
$\GL_\dd(\kk)$-stable subvariety $\MM$ 
of the representations, the space of
{\it nilpotent} representations).
We define:
$$
\MM =\MM^\dd:=\{ (M_1,\ldots,M_h)\mid 
M_1M_2\cdots M_h\in M_{d(h)\times d(h)}(\kk)
\ \mbox{is nilpotent} \}.
$$
The condition is equivalent to 
$M_{j\sh+1}M_{j\sh+2}\cdots M_h 
M_1 \cdots M_j\in M_{d(j)\times d(j)}(\kk)$ 
being nilpotent for any $j$.
In general $\MM$ is a connected but reducible variety.
I believe it has at most $h$ components, all of
equal dimension.
\\[1em]
{\bf Examples.}  (i) For $\dd=(1,1,1)$,
we have $\MM=\{(m_1,m_2,m_3)\in\kk^3\mid m_1m_2m_3=0\}$,
the union of the three coordinate planes.
Similarly for $\dd=(1^n)$.\\
(ii) For $\dd=(2,1,1)$, we have 
$\MM=\{([m_1,m_1'],[m_2,m_2']^T\!,m_3)\mid
(m_2m_2+m_1'm_2')\,m_3=0\}$, with two irreducible components
of dimension four.
\\[1em]
We define an isomorphism from $\MM$ 
to a union of opposite cells of Schubert
varieties in $\Fl(\dd,V)$.  Here we take
$n:=d_1+\cdots+d_h,$ so that $\dd$ is 
a composition of $n$, and we consider
$V=V_1\oplus\cdots\oplus V_n$, where $V_j=F^{d(j)}$. 

First, we embed 
$$
M_{d(h)\times d(1)}(\kk)\times\cdots
\times M_{d(h\sh-1)\times d(h)}(\kk)
\,\hookrightarrow\, G=GL_n(F),
$$
$$
M=(M_1,\ldots,M_h)\longmapsto \tM:=
\left(\begin{array}{@{\!}c@{\!}c@{\!}c@{\!}c@{\!}c@{\!}}
\, 0 & M_2 &\, 0 &\,\cdots &\, 0\,\\
\,0 & 0 &\, M_3 &\,\cdots &\, 0\,\\[-.5em]
\,\cdot &\cdot &\,\cdot&\, &\,\cdot\,\\[-.7em]
\,\cdot &\cdot &\,\cdot&\,&\,\cdot\,\\[-.3em]
\,0 & 0 &\, 0 &\,\cdots &\, M_h\, \\
\,t^{\sh-1}\!M_1 & 0 &\, 0 &\,\cdots &\,0
\end{array}\right).
$$
That is, $\tM=\tM_1+\cdots+\tM_h$, 
where $\tM_i:V\to V$,
$$
\tM_j(v):=\left\{\begin{array}{cl}
t^{-\delta_{1j}}M_j(v) &\ \mbox{for}\, v\in V_j\\[.2em]
0 &\ \mbox{for}\, v\in V_k,\, k\neq j.
\end{array}\right.
$$
We adopt the notations $\tM_{j+hk}:=\tM_j$, and:
$$
\tM_j^{[k]}:=\underbrace{
\tM_{j\sh-k\sh+1}\cdots\tM_{j\sh-1}\tM_j}_
{\mbox{\footnotesize $k$ factors}}\, .
$$
Note that $\tM_j\tM_k=0$ unless $k=j\sh+1$,
and for $M\in\MM$, we have
$\tM_j^{[h d(j)]}=(\tM_j^{[h]})^{d(j)}=0$,
so that $\tM^n=0$.
Now we can define $\psi_M:V\to V$,
\\[-.5em]
$$
\begin{array}{rcl}
\psi_M(v)&\!\!\!:=\!\!\!&
\displaystyle\frac{t^{n\sh-1}}{1-\tM}(v)\\[1em]
&\!\!\!=\!\!\!& 
t^{n\sh-1} (\ v_1+\,\tM_1(v_1)+
\tM_h\tM_1(v_1)+\cdots+\tM_1^{[nh\sh-h]}(v_1)\\[.3em]
&&\mbox{}\qquad +\, v_2+\,\tM_2(v_2)+
\tM_1\tM_2(v_2)+\cdots+\tM_2^{[nh\sh-h\sh+1]}(v_2)\\[.3em]
&&\mbox{}\qquad +\quad \cdots \\[.3em]
&&\mbox{}\qquad +\, v_h+\,\tM_h(v_h)+
\tM_{h\sh-1}\tM_h(v_h)+\cdots+\tM_h^{[nh\sh-1]}(v_h)\ )\ ,
\end{array}
$$
where $v=v_1+\cdots+v_h$ with $v_j\in V_j$.

Recall $E_{(\bullet)}=
(E_{(1)}\sh\supset\cdots\sh\supset E_{(h)})$,
the standard flag in $\Fl(\dd,V)$, where
$E_{(j)}:=E_{1+d(1)+\cdots+d(j\sh-1)}$ and 
$E^{(j)}\oplus E'_{(j)}=V$.
Then Lusztig's isomorphism is given by the map:
$$
\begin{array}{cccc}
\Psi:&\MM&\to&\Fl(\dd,V)\\[.1em]
&M&\mapsto&
\psi_M(E_{(\bullet)})
\end{array}
$$
where
$$
\psi_M(E_{(\bullet)}):=
(\,\psi_M(E_{(1)})\supset\cdots\supset\psi_M(E_{(h)})\,)\ .
$$

We give three coordinate descriptions of $\Psi$.
First, consider the decomposition,
$E_1=\kk^{d(1)}\oplus\cdots\oplus\kk^{d(h)}
\oplus t\kk^{d(1)}\oplus\cdots$,
so that we can write
$V\ni v=u_1+\cdots+u_h+tu_{h\sh+1}+\cdots$
with $u_i\in\kk^{d(i)\mod n}$.
Then we may write $\Psi(M)=
(\Lam_1\sh\supset\cdots\sh\supset\Lam_h)$
with
$$
\Lam_j=\{
u_1\sh+\cdots\sh+u_h\sh+tu_{h\sh+1}\sh+\cdots\ \mid\
u_{i\sh-1}= M_i(u_i)\ \ \forall\,i
\leq nh\sh-h\sh+j 
\}\, .
$$

Second, we write a partial flag 
$(\Lam_1\sh\supset\cdots\sh\supset\Lam_h)$
by a semi-infinite matrix of $n$ 
column-vectors $[v_1,\cdots,v_n]$ which are compatible
with all the lattices in the flag: that is,
$\Lam_j=\Span_\kk\langle v_i\rangle_
{i\geq 1+d(1)+\cdots+d(j\sh-1)}$.
We will write $[v_1,\cdots,v_n]$ as a block
matrix with blocks of sizes $d_1,\cdots,d_h$.
Let $I_j$ be an identity matrix of size $d_j$,
and denote $M_{j+hk}:=M_j$,\ \
$M_j^{[k]}:=M_{j\sh-k\sh+1}\cdots M_{j\sh-1}M_j$.
Then: 
$$
\Psi(M)\quad=\quad
\left\lceil\,\begin{array}{@{\!}c@{\!}c@{\!}c@{\!}c@{\!}}
M_1^{[nh\sh-h]}&\,M_2^{[nh\sh-h\sh+1]}&\,\cdots&
\ M_h^{[nh\sh-1]}\\[-.5em] 
\vdots&\vdots&&\vdots\\
M_1&\,M_1M_2&\,\cdots& M_h^{[h]}\\[.3em]
I_1&\,M_2&\,\cdots&\ M_h^{[h\sh-1]}\\[.3em]
0&\,I_2&\,\cdots&\ M_h^{[h\sh-2]}\\[-.4em] 
\vdots&\vdots&&\vdots\\
0&\,0&\,\cdots&\,I_h\\
0&\,0&\,\cdots&\,0\\[-.4em]
\vdots&\vdots&&\vdots\\
\end{array}\,\right\rceil
$$

Third, taking bases adapted to each lattice $\Lam_j$
(i.e. performing column reduction on the above matrix),
we obtain:
$$
\Lam_1=
\left\lceil\,\begin{array}{@{\!}c@{\!}c@{\!}c@{\!}c@{\!}}
M_1^{[nh-h]}&\ 0&\ \cdots&\ 0\\[-.4em]
\vdots&\vdots&&\vdots\\
M_hM_1&\ 0&\ \cdots&\ 0\\
M_1&\ 0&\ \cdots&\ 0\\
I_1&\ 0&\ \cdots&\ 0\\
0&\ I_2&\ \cdots&\ 0\\[-.4em]
\vdots&\vdots&&\vdots\\
0&\ 0&\ \cdots&\ I_h\\
0&\ 0&\ \cdots&\ 0\\[-.4em]
\vdots&\vdots&&\vdots\\&&&
\end{array}\ \right\rceil\, ,
\qquad
\Lam_2=
\left\lceil\,\begin{array}{@{\!}c@{\!}c@{\!}c@{\!}c@{\!}}
M_2^{[nh\sh-h\sh+1]}&\ 0&\ \cdots&\ 0\\[-.4em]
\vdots&\vdots&&\vdots\\
M_hM_1M_2&\ 0&\ \cdots&\ 0\\
M_1M_2&\ 0&\ \cdots&\ 0\\
M_2&\ 0&\ \cdots&\ 0\\
I_2&\ 0&\ \cdots&\ 0\\
0&\ I_3&\ \cdots&\ 0\\[-.4em]
\vdots&\vdots&&\vdots\\
0&\ 0&\ \cdots&\ I_1\\
0&\ 0&\ \cdots&\ 0\\[-.4em]
\vdots&\vdots&&\vdots
\end{array}\ \right\rceil\, ,
\qquad \mbox{etc.}
$$

From all these descriptions, we see that
$\Psi$ is $\GL_\dd(\kk)$-equivariant
(provided we diagonally embed
$\GL_\dd(\kk)\subset \GL_n(\kk)$).  
Also $E_1\supset \psi_M(E_{(1)})\supset\cdots
\supset\psi_M(E_{(h)})
\supset t^n E_1$ for all $M\in\MM$.

\subsection{Image of Lusztig's isomorphism}

The map $\Psi$ embeds $\MM$ in $\Fl(\dd,V)$.
We give several descriptions of the image.
From the last coordinate description of $\Psi$,
we may easily show:
$$
\Psi(\MM)=
\left\{\Ldot\in\Fl(\dd,V)
\left|\begin{array}{c}
\begin{array}{c@{\!}c@{\!}c@{\!}c@{\!}c@{\!}c@{\!}c@{\!}c@{\!}c}
\Lam_1&\stackrel{d(1)}{\supset}&\Lam_2&
\stackrel{d(2)}{\supset}&\ \, \cdots&\!\stackrel{d(h\sh-1)}{\supset}&
\Lam_h&\stackrel{d(h)}{\supset}&\ t\Lam_1\\
\mbox{\tiny $d(1)$}\,\cup\quad &&
\mbox{\tiny $d(2)$}\,\cup\quad &&&&
\mbox{\tiny $d(h)$}\,\cup\quad \\[.2em]
t^{n\sh-1}\!E_{(2)}\ &\supset&\ t^{n\sh-1}\!E_{(3)}\ 
&\supset&\ \, \cdots&\supset&t^n\!E_{(1)}
&=&t^{n\sh-1}\!E_{(h\sh+1)}
\end{array}\\[2.5em]
\Lam_1\cap t^{n\sh-1}\sh E'_{(2)}=\cdots
=\Lam_h\cap t^{n\sh-1}\sh E'_{(h\sh+1)}=0
\end{array}\right.\right\}\ .
$$

This suggests how to describe the image
as a union of opposite cells of Schubert varieties:
$\Psi(\MM)=\cup_{\pi} X'_{\pi}$ 
for certain $\pi\in\tW/W_\dd$.
Let $\ZZ_{(j)}:=\ZZ_{\geq 1+d(1)+\cdots+d(j\sh-1)}$,
and consider the sets
$\pi\ZZ_{(j)}$,\ $1\leq j\leq h$, which determine
$\pi$ modulo $W_\dd$.
These sets should contain numbers as small as possible
subject to the conditions:
$$
\begin{array}{c@{\!}c@{\!}c@{\!}c@{\!}c@{\!}c@{\!}c@{\!}c@{\!}c}
\pi\ZZ_{(1)}&\stackrel{d(1)}{\supset}&\pi\ZZ_{(2)}&
\stackrel{d(2)}{\supset}&\ \cdots&
\stackrel{d(h\sh-1)}{\supset}&
\pi\ZZ_{(h)}&\stackrel{d(h)}{\supset}&\ \tau\pi\ZZ_{(1)}\\
\mbox{\tiny $d(1)$}\,\cup\ &&
\mbox{\tiny $d(2)$}\,\cup\ &&&&
\mbox{\tiny $d(h)$}\,\cup\ \\[.2em]
\tau^{n\sh-1}\ZZ_{(2)}&\ \supset\ 
&\tau^{n\sh-1}\ZZ_{(3)}&\ \supset\ & 
\cdots&\supset&\tau^n\ZZ_{(1)}
\end{array}\ ,
$$
where $A\stackrel{d}{\supset} B$ means $A\supset B$
and $\#(A\sh\setminus B)=d$.
There should be at most $h$ such permutations $\pi$ which
are Bruhat-maximal.  One can construct them by first
maximizing a particular $\pi\ZZ_{(j)}$, 
then constructing the rest of the sets, which might
not always be possible.

%We will construct the sets $\pi\ZZ_{(j)}$:
%first, we define $\pi\ZZ_{(1)}$, 
%with the elements listed in an unusual order;
%then we remove $d_1$ of the elements, keeping
%the order of the rest; etc.
%Specifically, we let
%$$
%\begin{array}{rcl}
%\pi\ZZ_{(1)}:=&\{\hspace{-1.1em}&
%(n\sh-1)n\sh+d_1\sh+1,\, (n\sh-2)n\sh+d_1\sh+1,\,\ldots,\,
%(n\sh-d_1\sh-1)n\sh+d_1\sh+1,\\
%&&(n\sh-1)n\sh+d_1\sh+2,\, (n\sh-1)n\sh+d_1\sh+3,\, 
%(n\sh-1)n\sh+d_1\sh+4,\ldots
%\}\end{array}\, ,
%$$
%where the leftmost $d_1\sh+1$ elements are decreasing and 
%the rest increasing.  Next:
%$$
%\pi\ZZ_{(2)}:=
%\{\mbox{leftmost $d_2$ elements of $\pi\ZZ_{(1)}$}\}
%\ \cup\ \ZZ_{(3)}\, , 
%$$
%where ``leftmost'' means in the order listed above.
%The elements of $\pi\ZZ_{(2)}$ are again
%listed with some leftmost elements decreasing
%and the rest increasing.
%Now we continue in this manner:
%$$
%\pi\ZZ_{(j)}:=
%\{\mbox{leftmost $d_{j}$ 
%elements of $\pi\ZZ_{(j\sh-1)}$}\}
%\ \cup\ \ZZ_{(j\sh+1)}\, . 
%$$
%Finally, $\pi\!\!\mod W_\dd$ is completely defined by
%the blocks 
%$\pi\ZZ_{(j)}\setminus \pi\ZZ_{(j\sh+1)}.$

\vspace{1em}

\noindent
{\bf Example.} For $\dd=(1,1,1)$,\, $n=3$,
we have $\ZZ_{(j)}=\ZZ_{\geq j}$, and
the three irreducible components are
$\pi_1 = [5,9,10]=[2,3,1]\tau^{(1,2,3)}$;\
$\pi_2 = [8,6,10]=[2,3,1]\tau^{(2,3,1)}$;\
$\pi_3 = [9,8,7]=[2,3,1]\tau^{(2,2,2)}$.
Each $\pi_j$ is obtained by filling $\pi\ZZ_{(j)}$
with numbers as small as possible subject to 
$\pi\ZZ_{(j)}\stackrel{d(j)}{\supset} t^{n-\sh-1}E_{(j\sh+1)}$, 
then constructing the rest of the $\pi\ZZ_{(k)}$.
\\[-1.2em]

Specifically, to get $\pi=\pi_1$, we start with
$\pi\ZZ_{\geq 1}\stackrel{1}{\supset}\ZZ_{\geq 8}$,
yielding $\pi\ZZ_{\geq 1}=\{5,8,9,10,\ldots\}$.  Next 
$\pi\ZZ_{\geq 1}\stackrel{1}{\supset}
\pi\ZZ_{\geq 2}\stackrel{1}{\supset}\ZZ_{\geq 9}$,
yielding $\pi\ZZ_{\geq 2}=\{8,9,10,11,\ldots\}$.
Similarly $\pi\ZZ_{\geq 3}=\{8,10,11,12,\ldots\}$, 
and we conclude $\pi = [5,9,10]$.

The other components are obtained similarly.  Note that since
the $d_j$ are all equal to a constant $d=1$, 
the automorphism $\sigma^{d}=\sigma$ acts on $\Psi(\MM)$.  
In fact,  $\pi^{(1)}=\sigma \pi^{(2)} \sigma^{-1}
=\sigma^2 \pi^{(3)} \sigma^{-2}$.
\\[1em]
For a given $\dd$, the $GL_\dd(\kk)$-orbits of 
$\MM$ are distinguished from each other
by certain collections of invariants, the {\it rank numbers}
$\rr=(r_j^k)$, where 
$$
r_j^k:=\rank(M_j^{[k]}:
\kk^{d(j)}\to\kk^{d(j\sh-k)})
\quad\mbox{for}\quad 1\sh\leq j\sh\leq h,\ 
1\sh\leq k\sh\leq (n\sh-1)h\, .
$$
(In fact, it suffices to consider $k<h\sh\cdot\min(d_1,\ldots,d_h)$.)
We also define $r_j^0:=d_j$ and $r_j^k:=0$
if not otherwise defined.

For a given collection of rank numbers 
$\rr:=(r_{jk})_{jk}$, we denote the corresponding nilpotent
quiver orbit by $\MM^\rr\subset \MM$.  
Not every collection of rank numbers is realized by
an orbit:  rather, the nonnegative integer entries in 
$\rr$ must obey the constraints:
$$
m_j^k:=r_j^k-r_j^{k\sh+1}
-r_{j\sh-1}^{k\sh+1}+r_{j\sh-1}^{k\sh+2}\,\geq\, 0\, .
$$

[Note: $m_j^k$ is the multiplicity in $(M_1,\ldots,M_h)$
of the indecomposable quiver summand $I_j^k$ defined as follows:
letting $\bar i:=i\!\!\mod h$, we define vector spaces
$U_1,\ldots,U_h$ by $\oplus_{l=1}^h U_l:=
\Span_\kk\langle e_{j-i}\rangle_{0\leq i\leq k}$ with
$e_i\in U_{\bar i}$; and maps 
$L_{\bar i}\sh:U_{\bar i}\sh\to U_{\bar i\sh-1}$ with
$L_{\bar i}(e_i):=e_{i\sh-1}$ and 
$L_{\bar j-\bar k}(e_{j-k}):=0$.  
Thus  for $h=1$, the representation $I_j^k$ reduces to 
a nilpotent Jordan block of size $k\sh+1$.] 

The image of a quiver orbit under $\Psi$ can be described as:
\\[.5em]
$$
\Psi(\MM^\rr)=
\left\{\Ldot\in\Fl(\dd,V)
\left|\begin{array}{c}
\begin{array}{c@{\!}c@{\!}c@{\!}c@{\!}c@{\!}c@{\!}c@{\!}c@{\!}c}
\Lam_1&\stackrel{d(1)}{\supset}&\Lam_2&
\stackrel{d(2)}{\supset}&\ \, \cdots&\!\stackrel{d(h\sh-1)}{\supset}&
\Lam_h&\stackrel{d(h)}{\supset}&\ t\Lam_1\\
\hspace{-1.2em}\mbox{\tiny $d(1)$}\,\cup&&
\hspace{-1.2em}\mbox{\tiny $d(2)$}\,\cup&&&&
\hspace{-1.2em}\mbox{\tiny $d(h)$}\,\cup\\[.2em]
t^{n\sh-1}\!E_{(2)}\ &\supset&\ t^{n\sh-1}\!E_{(3)}\ 
&\supset&\ \, \cdots&\supset&t^n\!E_{(1)}
\end{array}\\[2.5em]
\Lam_1\cap E'_1=\cdots=\Lam_h\cap E'_h=0\\[.5em]
\dim(\Lam_j/\Lam_j\sh\cap E_{(k)})=r_j^{(nh-h+j-k+1)}\ \ 
\forall j,k
\end{array}\right.\right\}
$$
\vspace{.5em}

\noindent  Here we take $E_{(j+hk)}:=t^k E_{(j)}$.
The image of the orbit closure $\overline{\MM^\rr}$ is
obtained by replacing $=$ by $\leq$ in the Schubert
conditions.  
From this, we may deduce that 
$\Psi(\MM^\rr)=X^{\circ\prime}_\pi$ and
$\Psi(\overline{\MM^\rr})=X'_\pi$ 
for a certain explicitly constructable
$\pi=\pi^\rr\in\tW/W_\dd$.

\section{The Variety of circular complexes}

\subsection{Circular complexes and Lusztig's isomorphism}

We apply the previous constructions to a particularly simple,
but interesting case, intensively considered from a different 
point of view by Mehta and Trivedi \cite{mehta}.
For positive integers $a\leq b$, 
we consider the variety of {\it two-step circular complexes}
or {\it loop-complexes}:
$$
\LL =\LL_{a,b}:=\{ (X,Y)\in M_{b\times a}(\kk)\sh\times M_{a\times b}(\kk)
\mid XY=0,\ YX=0\}.
$$
Recall that any finite linear chain-complex can be ``rolled up'' into
such a two-step complex by letting $\kk^a$ (resp.~$\kk^b$) be
the direct sum of all the odd-numbered (resp.~even-numbered) spaces
in the linear complex.  This gives a natural map from the variety 
of chain-complexes to $\LL$.  

Now, $\LL$ is a subvariety of 
the representations of the affine quiver $\widehat A_1$;
a subvariety which is invariant under the 
natural action of the group 
$GL_{a,b}(\kk):=GL_a(\kk)\times GL_b(\kk)$, namely
$(g_a,g_b)\cdot (X,Y):=(g_bXg_a^{-1},g_aYg_b^{-1}).$
We easily see that $\LL$ is a finite union of 
$GL_{a,b}(\kk)$-orbits.  In fact, $\LL$
has exactly $a\sh+1$ open orbits $\LL^\circ_0,\ldots,\LL^\circ_a$, 
whose closures give the $a\sh+1$ irreducible components of $\LL$:
$$
\LL^\circ_c:=\{(X,Y)\in\LL\mid\, \rank X=c,\ \rank Y=a-c\},
\qquad \LL_c:=\overline{\LL^\circ_c}\, .
$$

We define an isomorphism from $\LL$ 
to a union of opposite cells of Schubert
varieties in the partial affine flag variety
$\Fl(a,b;V)$, where $V=F^n$ and $n=a+b$.  
Our notation will emphasize 
the block decomposition $V=F^a\oplus F^b$,
as well as:
$$
E=\kk^a\oplus \kk^b\oplus t\kk^a\oplus t\kk^b\oplus\cdots\ .
$$
In this case, Lusztig's isomorphism is given by the map
$\Psi:\LL\to G/P_a\cong \Fl(a,b;V)$ as:
$$
\Psi(X,Y):=
\left(\,\begin{array}{@{\!}c@{\!}c@{\!}}
tI_a&\ \ tY\\[.2em]
X&\ \ tI_b
\end{array}\,\right) \!\!\!\mod P_{\widehat{0,a}}
$$
where $I_m$ is an identity matrix of size $m$;
or in terms of lattices, $\Psi(X,Y)=(\Lam_1\supset\Lam_2)$,
where $E\supset\Lam_i$ and:
$$
\Lam_1\ =\
\left\lceil\begin{array}{@{\!}c@{\!}c@{\!}}
0&\ 0\\
X&\ 0\\
I_a&\ Y\\
0&\ I_b\\
0&\ 0\\[-.4em]
\vdots&\vdots\\ &
\end{array}\right\rceil
\ = \
\left\lceil\begin{array}{@{\!}c@{\!}c@{\!}}
0&\ 0\\
X&\ 0\\
I_a&\ 0\\
0&\ I_b\\
0&\ 0\\[-.4em]
\vdots&\vdots\\ &
\end{array}\right\rceil \!\!\mod P\ ,
\qquad
\Lam_2\ =\
\left\lceil\begin{array}{@{\!}c@{\!}c@{\!}}
0&\ 0\\
0&\ 0\\
Y&\ 0\\
I_b&\ X\\
0&\ I_a\\
0&\ 0\\[-.4em]
\vdots&\vdots
\end{array}\right\rceil
\ =\
\left\lceil\begin{array}{@{\!}c@{\!}c@{\!}}
0&\ 0\\
0&\ 0\\
Y&\ 0\\
I_b&\ 0\\
0&\ I_a\\
0&\ 0\\[-.4em]
\vdots&\vdots
\end{array}\right\rceil \!\!\mod P\ .
$$
Here the column vectors of each matrix
give an A-basis for $\Lam_i\subset E$, written with respect to
$E=\Span_\kk\langle e_i\rangle_{i\sh\geq 1}$.  The blocks
have sizes $a,b,a,b,\ldots$.  
The map $\Psi$ is $GL_{a,b}(\kk)$-equivariant,
provided we embed $GL_{a,b}(\kk)\subset GL_n(\kk)\subset G$
as block-diagonal matrices with constant coefficients.

We easily deduce:
$$
\Psi(\LL_c)=
\left\{\Ldot\in\Fl(a,b;V)
\left|\!\begin{array}{c}
\begin{array}{c@{\!}c@{\!}c@{\!}c@{\!}c@{\!}c}
E_{(2)}\ &\supset&\ E_{(3)}\ \\
\stackrel{b}{}\!\cup\quad &
&\stackrel{a}{}\!\cup\quad\\[-.6em]
\Lam_1&\csupset{a}&\Lam_2&\csupset{b}&\ \ t\Lam_1\\
\stackrel{a}{}\!\cup\quad &
&\stackrel{b}{}\!\cup\quad\\[.1em]
E_{(4)}\ &\supset&\ E_{(5)}\
\end{array}\ \
\begin{array}{c}
\Lam_1\sh\cap E'_{(3)}=\Lam_2\sh\cap E'_{(4)}=0\\[1em]
\dim(\Lam_1/\Lam_1\sh\cap E_{(3)})\leq c\\[.3em]
\dim(\Lam_2/\Lam_2\sh\cap E_{(4)})\leq a\sh-c
\end{array}
\end{array}\right.\!\!\!\!\!\!\right\}\ .
$$
Here 
$$
E_{(1)}\sh=E_1,\ \, 
E_{(2)}\sh=E_{a+1},\ \, 
E_{(3)}\sh=E_{a+b+1},\ \,  
E_{(4)}\sh=E_{2a+b+1},\ \,  
E_{(5)}\sh=E_{2a+2b+1}\  .
$$
From this, we may also realize $\LL_c$ 
as a subset of an ordinary flag variety 
$\Fl(b,a,b;\kk^{a+2b})$, the variety of 
partial flags 
$$
\kk^{a+2b}\csupset{b}U_1\csupset{a}U_2\csupset{b}0\ .
$$
In fact, let $\kk^{a+2b}=E_{(2)}/E_{(5)}$, with $\kk$-basis
$\{\bar e_{a\sh+1},\ldots, \bar e_{2a+2b}\}$,
where $\bar e_i:=e_i\!\!\mod E_{(5)}$.
Thus $\bar E_{(3)}\supset \bar E_{(4)}$ is the standard
flag in $\Fl(b,a,b;\kk^{a+2b})$.
We also define the nilpotent linear operator 
$\bar t:\kk^{a+2b}\to\kk^{a+2b}$ by 
$
\bar t(\bar e_i)=\bar e_{i\sh+n} \!\!\mod E_{(5)}\ .
$
Then we have:
$$
\LL_c\cong\Psi(\LL_c)\cong
\left\{
U_{\bullet}\in\Fl(b,a,b;\kk^{a+2b})
\left|\begin{array}{c}
U_1\supset \bar E_{(4)},\ \
U_2\subset \bar E_{(3)},\\[.3em]
\dim(U_1\sh\cap \bar E_{(3)})\geq a\sh+b\sh-c\\[.3em]
\dim(U_2\sh\cap \bar E_{(4)})\geq b\sh-a\sh+c\\[.3em]
U_1\sh\cap \bar E'_{(3)}=U_2\sh\cap \bar E'_{(4)}=0\\[.4em]
U_2\supset \bar t(U_1)
\end{array}
\right.\right\}
$$
This is precisely the opposite cell of a Schubert
variety in $\Fl(b,a,b;\kk^{a+2b})$, but with the
additional algebraic incidence condition 
$U_2\supset \bar t(U_1)$, which can be written
in terms of the Plucker coordinates
of $U_1, U_2$.

\subsection{Affine permutations for circular complexes}

We wish to identify the image of $\LL_c$ 
as the opposite cell of an affine Schubert variety, 
$$
\Psi(\LL_c)=X'_\pi\,, \quad\mbox{for some}
\quad\pi=\pi_c\in W_a\backslash\tW/W_a\,.
$$
First, we construct the sets
$\pi\ZZ_{(1)}, \pi\ZZ_{(2)}$, which then determine
$\pi$ modulo $W_a$.
These sets should contain numbers as small as possible
subject to the conditions:
$$
\begin{array}{c}
\begin{array}{c@{\!}c@{\!}c@{\!}c@{\!}c@{\!}c}
\ZZ_{(2)}\ &\supset&\ \ZZ_{(3)}\ \\
\stackrel{b}{}\!\cup\quad &
&\stackrel{a}{}\!\cup\quad\\[-.6em]
\pi\ZZ_{(1)}&\ \ \csupset{a}\ \  &\pi\ZZ_{(2)}&\ \ 
\csupset{b}\ \ &\tau\pi\ZZ_{(1)}\\
\stackrel{a}{}\!\cup\quad &
&\stackrel{b}{}\!\cup\quad\\[.1em]
\ZZ_{(4)}\ &\supset&\ \ZZ_{(5)}\
\end{array}\qquad
\begin{array}{c}
\#(\pi\ZZ_{(1)}\sh\setminus \ZZ_{(3)})=c\\[.3em]
\#(\pi\ZZ_{(2)}\sh\setminus \ZZ_{(4)})=a\sh-c
\end{array}\ ,
\end{array}
$$
To construct $\pi$ according to these constraints,
we will divide $[1,n]$ into intervals (blocks) of the form:  
$$
\underbrace{i,\ldots,j}_{(k)}:= [i,\,i\sh+1,\ldots,j],
$$
where $k=j\sh-i\sh+1$ is the number of integers in the
interval.  We perform two subdivisions as follows:
$$
\begin{array}{rll}
[1,\ldots,n]:= [\!\!\! &
\underbrace{1,\ldots,c}_{\mbox{\sc i}\ (c)},
\underbrace{c\sh+1,\ldots,a}_{\mbox{\sc ii}\ (a-c)},
\underbrace{a\sh+1,\ldots,2a\sh-c}_{\mbox{\sc iii}\ (a-c)},
\\[2em]&
\underbrace{2a\sh-c\sh+1,\ldots,a\sh+b\sh-c}_{\mbox{\sc iv}\ (b-a)},
\underbrace{a\sh+b\sh-c\sh+1,\ldots,a\sh+b}_{\mbox{\sc v}\ (c)}
\ \, ] \end{array}\ ,
$$
$$
\begin{array}{rll}
[1,\ldots,n]:= [\!\!\! &
\underbrace{1,\ldots,a\sh-c}_{\mbox{\sc i}'\ (a-c)},
\underbrace{a\sh-c\sh+1,\ldots,a}_{\mbox{\sc ii}'\ (c)},
\underbrace{a\sh+1,\ldots,a\sh+c}_{\mbox{\sc iii}'\ (c)},
\\[2em]&
\underbrace{a\sh+c\sh+1,\ldots,b\sh+c}_{\mbox{\sc iv}'\ (b-a)},
\underbrace{b\sh+c\sh+1,\ldots,a\sh+b}_{\mbox{\sc v}'\ (a-c)}
\ \, ] \end{array}\ ,
$$
where we have numbered the blocks with roman numerals.
Now $\pi$ takes the first set of blocks to
the second set, as well as shifting them
by powers of $\tau$:
$$
\begin{array}{rll}
\pi:= [\!\!\! &
\underbrace{a\sh+1,\ldots,a\sh+c}_{\mbox{\sc v}'\ (c)},
\underbrace{a\sh+2b\sh+c\sh+1,\ldots,2a\sh+2b}_{\tau(\mbox{\sc iii}')\ (a-c)},
\underbrace{a\sh+b\sh+1,\ldots,2a\sh+b\sh-c}_{\tau(\mbox{\sc ii}')\ (a-c)},
\\[2em]&
\underbrace{2a\sh+b\sh+c\sh+1,\ldots,a\sh+2b\sh+c}_{\tau(\mbox{\sc iv}')\ (b-a)},
\underbrace{3a\sh+2b\sh-c\sh+1,\ldots,3a\sh+2b}_{\tau^2(\mbox{\sc i}')\ (c)}
\ \, ] \end{array}\ .
$$
Recall that $\pi$ represents the double coset 
$W_a\pi\,W_a$: in fact, $\pi$ is maximal with respect to
the left action of $W_a$, and minimal with respect to
the right action.  This is the correct normalization
so that $\ell(\pi_c)=\dim_\kk(\LL_c)$.

To analyze the decomposition of $\pi$ into simple reflections,
we construct its loop wiring diagram.
As before, the strip below (with top and bottom edges identified)
represents a cylinder  with $n=a+b$ dots on either end.
For each $i$ we write $\pi(i)=\bar\pi(i)+nj$, 
and we draw a wire connecting the dot $i$
on the right to the dot $\bar\pi(i)$ on the left,
but looping upwards (around the cylinder) $j$ times.
We will group the wires into five {\it cables} corresponding
to our blocks {\sc i,\!\ldots,v} (on the right)
and $\mbox{\sc i}',\ldots,\mbox{\sc v}'$ (on the left),
so that the cable starting at {\sc i} represents $c$
non-crossing wires, etc.  As a final simplification,
instead of drawing the diagram for $\pi$, 
we instead draw the diagram for $\tau^{-1}\pi$
(a harmless normalization, since $\tau$ is in the center of
$\tW$).
\\[1.5em]
\mbox{}\hspace{.2in}
%TexCad Options
\unitlength .85mm
\linethickness{0.4pt}
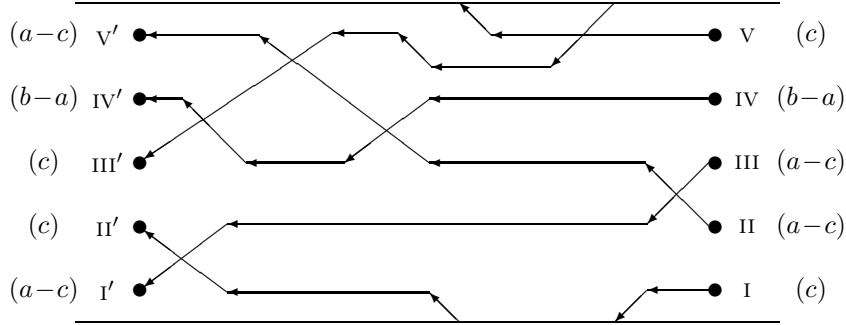
\begin{picture}(125.00,55.00)
\put(20.00,10.00){\circle*{2.11}}
\put(110.00,10.00){\circle*{2.11}}
\put(20.00,20.00){\circle*{2.11}}
\put(110.00,20.00){\circle*{2.11}}
\put(20.00,30.00){\circle*{2.11}}
\put(110.00,30.00){\circle*{2.11}}
\put(20.00,40.00){\circle*{2.11}}
\put(110.00,40.00){\circle*{2.11}}
\put(20.00,50.00){\circle*{2.11}}
\put(110.00,50.00){\circle*{2.11}}
\put(10.00,5.00){\line(1,0){110.00}}
\put(10.00,55.00){\line(1,0){110.00}}
\put(15.00,50.00){\makebox(0,0)[cc]{$\mbox{\sc v}'$}}
\put(15.00,40.00){\makebox(0,0)[cc]{$\mbox{\sc iv}'$}}
\put(15.00,30.00){\makebox(0,0)[cc]{$\mbox{\sc iii}'$}}
\put(15.00,20.00){\makebox(0,0)[cc]{$\mbox{\sc ii}'$}}
\put(15.00,10.00){\makebox(0,0)[cc]{$\mbox{\sc i}'$}}
\put(115.00,50.00){\makebox(0,0)[cc]{\mbox{\sc v}}}
\put(115.00,40.00){\makebox(0,0)[cc]{\mbox{\sc iv}}}
\put(115.00,30.00){\makebox(0,0)[cc]{\mbox{\sc iii}}}
\put(115.00,20.00){\makebox(0,0)[cc]{\mbox{\sc ii}}}
\put(115.00,10.00){\makebox(0,0)[cc]{\mbox{\sc i}}}
\put(5.00,50.00){\makebox(0,0)[cc]{$(a\!-\!c)$}}
\put(5.00,40.00){\makebox(0,0)[cc]{$(b\!-\!a)$}}
\put(5.00,30.00){\makebox(0,0)[cc]{$(c)$}}
\put(5.00,20.00){\makebox(0,0)[cc]{$(c)$}}
\put(5.00,10.00){\makebox(0,0)[cc]{$(a\!-\!c)$}}
\put(125.00,50.00){\makebox(0,0)[cc]{$(c)$}}
\put(125.00,40.00){\makebox(0,0)[cc]{$(b\!-\!a)$}}
\put(125.00,30.00){\makebox(0,0)[cc]{$(a\!-\!c)$}}
\put(125.00,20.00){\makebox(0,0)[cc]{$(a\!-\!c)$}}
\put(125.00,10.00){\makebox(0,0)[cc]{$(c)$}}
\put(109.00,30.00){\vector(-1,-1){9.67}}
\put(99.33,20.33){\vector(-1,0){65.67}}
\put(33.67,20.33){\vector(-4,-3){13.00}}
\put(109.00,20.00){\vector(-1,1){10.00}}
\put(99.00,30.00){\vector(-1,0){33.67}}
\put(65.33,30.00){\vector(-4,3){26.67}}
\put(38.67,50.00){\vector(-1,0){17.67}}
\put(109.00,40.00){\vector(-1,0){43.67}}
\put(65.33,40.00){\vector(-4,-3){13.33}}
\put(52.00,30.00){\vector(-1,0){15.33}}
\put(36.67,30.00){\vector(-1,1){10.00}}
\put(26.67,40.00){\vector(-1,0){5.67}}
\put(109.00,10.00){\vector(-1,0){9.67}}
\put(99.33,10.00){\vector(-1,-1){5.00}}
\put(94.33,55.00){\vector(-1,-1){10.00}}
\put(84.33,45.00){\vector(-1,0){18.67}}
\put(65.67,45.00){\vector(-1,1){5.33}}
\put(60.33,50.33){\vector(-1,0){10.00}}
\put(50.33,50.33){\vector(-3,-2){29.67}}
\put(109.00,50.00){\vector(-1,0){34.00}}
\put(75.00,50.00){\vector(-1,1){5.00}}
\put(70.00,5.00){\vector(-1,1){4.67}}
\put(65.33,9.67){\vector(-1,0){31.67}}
\put(33.67,9.67){\vector(-4,3){13.00}}
\end{picture}

\noindent Whenever a cable with $k$ wires crosses
one with $k'$ wires, we have a total of $kk'$ wire
crossings.  Thus the six cable-crossings of our picture
give a wire-crossing total of:
$$
\ell(\pi)=(a-c)^2+c^2+(a-c)(b-a)+2c(a-c)+c(b-a)=ab\, ,
$$
which we may confirm by checking directly that $\dim(\LL_c)=ab$.

Now we may write a reduced decomposition for $\pi$
as follows.  For integers $i,k$, define the affine permutation 
$s_i^{[k]}:=s_i s_{i-2}\cdots s_{i-2k+2}$, which has $k$
mutually commuting factors.  Recall our convention
$s_{i+nj}:=s_i$.  For each cable crossing:
$$
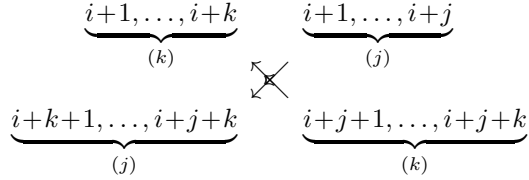
\begin{array}{r@{\!}c@{\!}c@{\!}c@{\!}l}
\underbrace{i\sh+1,\ldots,i\sh+k}_{(k)}&&&&
\underbrace{i\sh+1,\ldots,i\sh+j}_{(j)}\\[-.5em]
&\ \ \nwarrow&&\!\swarrow\ \ &\\[-.7em]
%&&\bullet&&\\[-.7em]
&\ \ \swarrow&&\!\nwarrow\ \ &\\[.1em]
\underbrace{i\sh+k\sh+1,\ldots,i\sh+j\sh+k}_{(j)}
&&&&
\underbrace{i\sh+j\sh+1,\ldots,i\sh+j\sh+k}_{(k)}
\end{array}
$$
we define the associated ``totally commutative'' 
permutation:
$$
s_{i\sh+1}^{[j,k]}\ :=\
s_{i\sh+k}^{[1]}s_{i\sh+k\sh+1}^{[2]}\cdots
s_{i\sh+k\sh+m}^{[\min(m,j,k)]}\cdots
s_{i\sh+j\sh+k}^{[\min(j,k)]}\cdots
s_{i\sh+j\sh+m'}^{[\min(m'\!,j,k)]}\cdots
s_{i\sh+j\sh+1}^{[2]}s_{i\sh+j}^{[1]}.
$$
Finally, we can write 
$$
\pi = 
\tau\,
s_{1}^{[a-c,c]}\,
s_{a+1}^{[b-a,c]}\,
s_{b+1}^{[a-c,c]}\,
s_{a+1}^{[a-c,b\sh-a]}\,
s_{b+a-c+1}^{[c,c]}\,
s_{c+1}^{[a-c,a-c]}\ ,
$$
where the six factors (other than $\tau$) correspond to 
the cable crossings, listed left to right.

\subsection{Bott-Samelson resolution}

We can use the above data to give a Bott-Samelson
resolution of singularities for $\LL_c$.  
Although this is clearly far from a minimal resolution,
it brings the circular complexes into the framework of 
Frobenius splittings and other results for
Bott-Samelson varieties (cf.~Mathieu \cite{mathieu},
Kumar \cite{kumar2}, Ramanathan \cite{ramanathan},\ldots).  
In particular, we have the following results proved
by Mehta-Trivedi \cite{mehta}
\\[1em]
{\bf Theorem}  {\it The variety of circular complexes
$\LL_c$ and the closures of all its $GL_n(\kk)$-orbits
are normal, Cohen-Macaulay, and have rational
singularities.}
\\[1em]
The construction of the affine Bott-Samelson variety
$Z_{\ii}$ corresponding to a reduced word 
$\ii=(i_1,\ldots,i_l)$
is exactly analogous to (and includes
as a special case) the construction for $GL_n(\kk)$.
(Cf.~\S1.4.)

We illustrate with the simplest example in our case: 
$c=1$, $a=2$, $b=3$, $n=5$ so that each of
the blocks $\mbox{\sc i},\ldots,\mbox{\sc v}$ has
size $1$, and our cable diagram in the previous
section is a simple wiring diagram.
Then $\pi=s_1 s_3 s_4 s_3 s_0 s_2$,
and the Bott-Samelson variety is:
$$
Z_{\ii}\ :=\
\left\{\begin{array}{c}
(\Lam_1,\Lam_2,\Lam_3,\Lam_4,\Lam_5,\Lam'_4)\\
\in \Gr(V)^6
\end{array}
\left|\begin{array}{c@{\!}c@{\!}c@{\!}c@{\!}c@{\!}
c@{\!}c@{\!}c@{\!}c@{\!}c@{\!}c}
&&&&&&E_4&&&\\
&&&&&\swarrow&&\nwarrow\\
E_1&\ \leftarrow\ &E_2&\ \leftarrow\ &E_3&\ \leftarrow\ 
&\Lam'_4&\ \leftarrow\ &E_5&\ \leftarrow\ &tE_1\\
&\nwarrow&&\swarrow&&\nwarrow&&\nwarrow&&\swarrow&\\
\Lam_1&\ \leftarrow\ &\Lam_2&\ \leftarrow\ &\Lam_3&\ \leftarrow\ &
\Lam_4&\ \leftarrow\ &\Lam_5&\ \leftarrow\ &t\Lam_1
\end{array}\right.\right\}.
$$
Here $E_j:=\Span_\kk\langle e_i\rangle_{i\geq j}$,
and each arrow $U\sh\leftarrow V$ indicates
the conditions $U\supset V$,\ $\dim_\kk(U/V)=1$.
We construct a point of $Z_\ii$ by starting with 
the standard flag $E_1\supset E_3\supset\cdots$,
and successively choosing the spaces $\Lam_2$,
$\Lam'_4$, $\Lam_5$, $\Lam_4$, $\Lam_1$, $\Lam_3$,
corresponding to the letters 1,3,4,3,0,2.
Each such choice corresponds to
a fibration with fiber $\PP^1$, hence $Z_\ii$
is smooth of dimension $\ell(\pi)=6$.  As we did for
$X_\pi$, we can embed
$Z_\ii$ into a finite-dimensional flag variety
for the $\kk$-vector space $t^{-1}\!E_3/tE_1$,
since $t^{-1}E_3\supset\Lam\supset tE_1$
for $\Lam=\Lam_1,\ldots,\Lam_5,\Lam'_4$.

We can define
a regular, birational map of $Z_\ii$
onto $X_\pi$ by
forgetting all the spaces except 
$(\Lam_1,\Lam_3)$.  This map is generically
one-to-one because generically all the spaces
are determined by $\Lam_1$, $\Lam_3$: that is,
$\Lam_2=\Lam_1\sh\cap E_1$,\ $\Lam_5=t\Lam_1\sh+tE_1$, etc.
To desingularize the opposite cell $X'_\pi\cong\LL_c$, 
we consider the subset of $Z_\ii$ where
$\Lam_1, \Lam_3$ are generic with 
respect to the opposite standard flag $E'_1, E'_3$.

\end{document}